# Differentiation in Bundles with a Hyperspace Base

**Mark Burgin**
Department of Mathematics, UCLA
405 Hilgard Ave.
Los Angeles, 90046, USA

**Abstract**: It is possible to perform some operations with extrafunctions applying these operations separately to each coordinate. Operations performed in this manner are called regular. It is proved that it is possible to extend several operations with functions to regular operations with extrafunctions. Examples of such operations are addition of real functions and multiplication of real functions by numbers. However, there are operations with functions the extension of which by coordinates does not work because their application is not invariant with respect to representations of extrafunctions. One of such operations is differentiation, which is important for calculus, differential equations and many applications.

In this work, a method of regularization of irregular operations is developed and applied to differentiation. The main constructions are put together in the context of fiber bundles over hyperspaces of differential vector spaces and differential algebras.

*Keywords*: vector bundle, section, differentiation, extrafunction, operation, regularization

## 1. Introduction

Extension of functions to extrafunctions is only the first step in building the theory of extrafunctions. One of the basic components of such a theory explains how to perform operations with extrafunctions and explores what properties these operations have. In many cases, it is possible to take an operation with ordinary functions and to perform similar operations with extrafunctions applying these function operations separately to each coordinate. Operations performed in this manner are called *regular*. Examples of such operations are addition of extrafunctions and multiplication of extrafunctions by real or complex numbers. For instance, if we have two extrafunctions $F = \text{EF}_Q(f_i)_{i \in \omega}$ and $G = \text{EF}_Q(g_i)_{i \in \omega}$, then their sum is defined as $F + G = \text{EF}_Q(f_i + g_i)_{i \in \omega}$, i.e., addition is performed separately in each coordinate $i$. It is important that, as it is proved in (Burgin, 2001), the result does not depend on the choice of the extrafunction representations with which addition is performed and thus, it correctly defines operations with extrafunctions.

However, there are operations with ordinary functions the extension of which by coordinates does not work because their application is not invariant with respect to the representation of extrafunctions. Such operations are called *irregular*. For instance, multiplication of extrafunctions is an irregular operation (Burgin, 2002; 2004). One of irregular operations is differentiation, which is important for calculus, differential equations and many applications.

In this work, a general method of regularization of irregular operations is developed and applied to differentiation. The main constructions are put together in the context of bundles over hyperspaces of differential vector spaces and differential algebras. This allows regularization of irregular operations for different classes of norm-based extrafunctions, which include distributions, hyperdistributions, restricted pointwise extrafunctions, and compactwise extrafunctions.

## 2. Norm-based Sequential Hyperspaces

Let us consider a linear space $L$ over the field $F$ with a set $Q$ of seminorms in $L$ and the set $L^\omega$ is the set of all sequences with elements from $L$. We denote sequences by $(f_i)_{i \in \omega}$ or by $\{ f_i; i \in \omega \}$ or by $\{ f_i; i = 1, 2, 3, \ldots \}$.

**Proposition 2.1.** $L^\omega$ is a linear space over the field $F$.

Indeed, for arbitrary sequences $f = (f_i)_{i \in \omega}$ and $g = (g_i)_{i \in \omega}$ from the set $L^\omega$, we define operations of addition and multiplication by elements from $F$:

$$f + g = (f_i + g_i)_{i \in \omega}$$

and if $a \in F$, then

$$cf = (cf_i)_{i \in \omega}$$

Besides, we can see that any identity that is valid in $L$ and involves only operations of addition and/or multiplication by elements from $F$ is also valid in $L^\omega$.

There is a natural injection $\mu_L: L \to L^\omega$ where for any element $f$ from $L$, $\mu_L(f) = (f_i)_{i \in \omega}$ and $f_i = f$ for all $i = 1, 2, 3, \ldots$. Such elements $(f_i)_{i \in \omega}$ are called *stable* (in $L^\omega$).

Definitions imply the following result.

**Proposition 2.2.** $\mu$ is a monomorphism of linear spaces over the field $F$.

To build the superspace from the space $L$, we introduce an equivalence relation.

**Definition 2.1.** For arbitrary sequences $f = (f_i)_{i \in \omega}$ and $g = (g_i)_{i \in \omega}$ from the set $L^\omega$,

$$f \approx_Q g \text{ means that } \lim_{i \to \infty} q(f_i - g_i) = 0 \text{ for any } q \in Q$$

**Lemma 2.1.** The relation $\approx_Q$ is an equivalence relation in $L^\omega$.

<u>Proof.</u> By definition, this relation is reflexive. Besides, it is symmetric because $q(f - g) = q(g - f)$ for any seminorm $q$ and we need only to show that the relation $\approx_Q$ is transitive. Taking three sequences $f = (f_i)_{i \in \omega}$, $g = (g_i)_{i \in \omega}$ and $h = (h_i)_{i \in \omega}$ of elements from $F$ such that $f \approx_Q g$ and $g \approx_Q h$, for any seminorm $q$ from $Q$, we have

$$\lim_{i \to \infty} q(f_i - g_i) = 0$$

and

$$\lim_{i \to \infty} q(g_i - h_i) = 0$$

By properties of seminorms and limits, we have

$$0 \le \lim_{i \to \infty} q(f_i - h_i) = \lim_{i \to \infty} q(f_i - g_i + g_i - h_i) \le$$
$$\lim_{i \to \infty} q(f_i - g_i) + \lim_{i \to \infty} q(g_i - h_i) = 0 + 0 = 0$$

Consequently,

$$\lim_{i \to \infty} q(f_i - h_i) = 0$$

i.e., $f \approx_Q h$.

Lemma is proved.

Classes of the equivalence relation $\approx_Q$ form the hyperspace $L_{\omega Q}$.

In such a way, any sequence $f = (f_i)_{i \in \omega}$ with elements from $L$ determines (and represents) an element $F = \mathrm{Hs}_Q(f_i)_{i \in \omega}$ from the hyperspace $L_{\omega Q}$ and this sequence is called a *defining sequence* or *representation sequence* or *representation* of the element $F = \mathrm{Hs}_Q(f_i)_{i \in \omega}$. There is also a natural projection $\pi_Q: L^\omega \to L_{\omega Q}$ where $\pi_Q((f_i)_{i \in \omega}) = \mathrm{Hs}_Q(f_i)_{i \in \omega}$ for any sequence $f = (f_i)_{i \in \omega} \in L^\omega$. An element $F$ from $L_{\omega Q}$ is called *stable* if it has a stable representation, i.e., $F = \mathrm{Hs}_Q(f_i)_{i \in \omega}$ for some stable element $f = (f_i)_{i \in \omega}$.

Let us consider some examples.

**Example 2.1**. Taking the space $R$ of all real numbers, $Q = \{\text{absolute value}\}$ and applying the construction of the hyperspace, we obtain the set of real hypernumbers $R_\omega$ (Burgin, 2004).

**Example 2.2**. Taking the space $C$ of all complex numbers, $Q = \{\text{absolute value}\}$ and applying the construction of the hyperspace, we obtain the set of complex hypernumbers $C_\omega$ (Burgin, 2005).

**Example 2.3**. It is possible to define seminorm $q_{\mathrm{pt}x}$ in the space $F(R, R)$ of all real functions using an arbitrary real number $x$ and the following formula

$$q_{\mathrm{pt}x}(f) = |f(x)|$$

The set $Q_{\mathrm{pt}} = \{q_{\mathrm{pt}x}; x \in R\}$ of seminorms determines the equivalence relation $\approx_{\mathrm{pt}}$ in the space of all sequences of real functions. Then the hyperspace $F(R, R)_{\omega Q_{\mathrm{pt}}}$ is isomorphic to the set $F(R, R_\omega)$ of all restricted real pointwise extrafunctions (Burgin, 2004).

**Example 2.4**. It is possible to define seminorm $q_{\mathrm{pt}x}$ in the space $F(R, C)$ of all functions from $R$ to $C$ using an arbitrary real number $x$ and the following formula

$$q_{\mathrm{pt}x}(f) = |f(x)|$$

The set $Q_{\mathrm{pt}} = \{q_{\mathrm{pt}x}; x \in R\}$ of seminorms determines the equivalence relation $\approx_{\mathrm{pt}}$ in the space of all sequences of real functions. Then the hyperspace $F(R, C)_{\omega Q_{\mathrm{pt}}}$ is isomorphic to the set $F(R, C_\omega)$ of all restricted complex valued pointwise extrafunctions.

**Example 2.5**. It is possible to define seminorm $q_{\mathrm{pt}x}$ in the space $F(C, C)$ of all complex functions using an arbitrary real number $x$ and the following formula

$$q_{\mathrm{pt}x}(f) = |f(x)|$$

The set $Q_{pt} = \{q_{ptx}; x \in C\}$ of seminorms determines the equivalence relation $\approx_{pt}$ in the space of all sequences of real functions. Then the hyperspace $F(C, C)_{\omega Q_{pt}}$ is isomorphic to the set $F(C, C_\omega)$ of all restricted complex pointwise extrafunctions.

**Example 2.6**. Let us consider a set **K** of test functions (cf., for example, (Rudin, 1973)) such that for any real function $f(x)$ and any function $g(x)$ from **K**, the integral $\int f(x)g(x)dx$ exists. Then it is possible to define the seminorm $q_g$ in **F** using an arbitrary function $g(x)$ from **K** and the following formula

$$q_g(f) = |\int f(x)g(x)dx|$$

The set $Q_{\mathbf{K}} = \{q_g; g \in \mathbf{K}\}$ of seminorms determines the equivalence relation $\approx_{\mathbf{K}}$ in the space of all sequences of real functions. Then taking the hyperspace $F(R, R)_{\omega Q_{\mathbf{K}}}$ with respect to the corresponding integration, we get real extended distributions studied in (Burgin, 2004), as well as real measure-wise extrafunctions introduced and studied in (Burgin 2002).

**Example 2.7**. It is possible to define seminorm $q_{max[a, b]}$ in the space $C(R, R)$ of all continuous real functions using an arbitrary interval $[a,b]$ and the following formula

$$q_{max[a, b]}(f) = \max \{|f(x)|; x \in [a, b]\}$$

The set $Q_{comp} = \{q_{max[a, b]}; a, b \in R\}$ of seminorms determines the equivalence relation $\approx_{comp}$ in the space of all sequences of continuous real functions. Then the hyperspace $C(R, R)_{\omega Q_{comp}}$ is to the set $Comp(R, R_\omega)$ of all real compactwise extrafunctions (Burgin, 2004).

**Example 2.8**. If **F** is a normed field, then it is possible to treat **F** as a one-dimensional vector space over itself. The previous constructions give us the space $F^\omega$ and the hyperspace $F_\omega$.

Let us consider operations in the hyperspace $L_{\omega Q}$.

**Proposition 2.3.** Operations of addition and/or subtraction are correctly defined in the hyperspace $L_{\omega Q}$.

<u>Proof</u>. Let us take two elements $F = Hs_Q(f_i)_{i \in \omega}$ and $G = Hs_Q(g_i)_{i \in \omega}$ from the hyperspace $L_{\omega Q}$. We define $F + G = Hs_Q(f_i + g_i)_{i \in \omega}$. To show that this is an operation in the hyperspace, it is necessary to prove that $F + G$ does not depend on the choice of a representing sequences $(f_i)_{i \in \omega}$ and $(g_i)_{i \in \omega}$ for the elements $F$ and $G$. To do this, let us take another sequence $(h_i)_{i \in \omega}$ that represents $F$ and show that $F + G = \mathrm{Hs}_Q(h_i + g_i)_{i \in \omega}$. Note that all three sequences $(f_i)_{i \in \omega}$, $(h_i)_{i \in \omega}$ and $(g_i)_{i \in \omega}$ belong to $\boldsymbol{L}^\omega$.

In this case, for any seminorm $q$ from $Q$, we have

$$\lim_{i \to \infty} q((f_i + g_i) - (h_i + g_i)) = 0$$

because $\lim_{i \to \infty} q(f_i - h_i) = 0$. Thus, addition is correctly defined for all elements from the hyperspace $\boldsymbol{L}_{\omega Q}$.

The proof for the difference of elements from the hyperspace $\boldsymbol{L}_{\omega Q}$ is similar.

Proposition is proved.

The construction of addition and subtraction of elements from the hyperspace $\boldsymbol{L}_{\omega Q}$ implies the following result.

**Proposition 2.4.** Any identity that involves only operations of addition and/or subtraction and is valid for elements from the space $\boldsymbol{L}$ is also valid for elements from the hyperspace $\boldsymbol{L}_{\omega Q}$.

Let us assume that the class **F** is closed with respect to multiplication by real numbers, i.e., if $a$ is a real number and $f \in \mathbf{F}$, then the product $af$ also belongs to **F**. Then this operation induces the corresponding operation in sets of represented in **F** extrafunctions.

**Proposition 2.5.** Multiplication by elements from the field $\boldsymbol{F}$ is correctly defined for elements from the hyperspace $\boldsymbol{L}_{\omega Q}$.

<u>Proof.</u> Let us take an element $F = \mathrm{Hs}_Q(f_i)_{i \in \omega}$ from the hyperspace $\boldsymbol{L}_{\omega Q}$ and a real number $a$. We define $aF = \mathrm{Hs}_Q(af_i)_{i \in \omega}$. Let us take another sequence $(h_i)_{i \in \omega}$ that represents $F$ and show that $aF = \mathrm{Hs}_Q(ah_i)_{i \in \omega}$. Indeed, for any seminorm $q$ from $Q$, we have

$$\lim_{i \to \infty} q((af_i - ah_i)) = \lim_{i \to \infty} |a| q((f_i - h_i)) = |a| \cdot \lim_{i \to \infty} q((f_i - h_i)) = 0$$

because $\lim_{i \to \infty} q(f_i - h_i) = 0$ and $q$ satisfies Condition N2. Thus, multiplication by elements from the field $\boldsymbol{F}$ is correctly defined in the hyperspace $\boldsymbol{L}_{\omega Q}$.

Proposition is proved.

Definitions imply the following result.

**Proposition 2.6.** Any identity that is valid for elements from the space $\boldsymbol{L}$ and involves only operations of addition and/or multiplication by real numbers is also valid for elements from the hyperspace $\boldsymbol{L}_{\omega Q}$.

In particular, Proposition 2.6 gives the following identities for extrafunctions.

(1) $F(x) + G(x) = G(x) + F(x)$

(2) $F(x) + (G(x) + H(x)) = (F(x) + G(x)) + H(x)$

(3) $a(F(x) + G(x)) = aF(x) + aG(x)$

(4) $(a + b)F(x) = aF(x) + bF(x)$

(5) $a(bF(x)) = (ab)F(x)$

(6) $1 \cdot F(x) = F(x)$

Propositions 2.1, 2.2 and 2.4 imply the following results.

**Theorem 2.1.** $L_{\omega Q}$ is a linear space over the field $F$.

**Theorem 2.2.** $\pi_Q$ is an epimorphism of the linear space $L^\omega$ onto the hyperspace $L_{\omega Q}$.

An important question is when the space $L^\omega$ is a natural subspace of the hyperspace $L_{\omega Q}$.

**Definition 2.2.** A system $Q$ of seminorms in $L$ *separates* $L$ if for any elements $f$ and $g$ from $L$, there is a seminorm $q$ from $Q$ such that $q(f) \neq q(g)$.

**Lemma 2.2.** Any norm in $L$ separates $L$.

It is possible to define a mapping $\mu_Q: L \to L_{\omega Q}$ by the formula $\mu_Q(f) = \pi_Q(\mu(f))$.

**Proposition 2.7.** $\mu_Q$ is a monomorphism (inclusion) of linear spaces over the field $F$ if and only if the system $Q$ of seminorms separates $L$.

Proof. Sufficiency. If $f$ and $g$ are elements from $L$ and the system $Q$ of seminorms in $L$ separates $L$, then there is a seminorm $q$ from $Q$ such that $q(f) \neq q(g)$. Consequently, by properties of seminorms, $q(f - g) \neq 0$. Then by definition, $\mu_Q(f) \neq \mu_Q(f)$. So, by Theorem 2.2 and Proposition 2.2, $\mu_Q$ is a monomorphism.

Necessity. If the system $Q$ of seminorms does not separate $L$, then for some elements $f$ and $g$ from $L$, $q(f) = q(g)$ for all seminorms $q$ from $Q$. Consequently, by properties of seminorms, $q(f - g) = 0$. It means, by Definition 2.1, that $\mu(f) \approx_Q \mu(g)$. Thus, $\mu_Q(f) = \pi_Q(\mu(f)) = \pi_Q(\mu(fg)) = \mu_Q(g)$. So, $\mu_Q$ is not a monomorphism (inclusion). Thus, by the Law of Contraposition for propositions (cf., for example, (Church, 1956)), the condition from Proposition 2.2 is necessary.

Proposition is proved.

The set $Q = \{q_t ; t \in K\}$ of seminorms in the space $L$ allows us to define a topology in the space $L$.

**Definition 2.3.** If $l$ is an element from $L$ and $N = \{r_t ; t \in K\}$ is a set of positive real numbers $r_t$, then a set of the form

$$O_N l = \{g \in L \,;\, \forall t \in K \,(\, q_t(l - g) < r_t)\}$$

is called a *locally defined neighborhood* of the element $l$.

**Lemma 2.3.** The system of all locally defined neighborhoods determines a topology $\upsilon_{lQ}$ called *local Q-topology* in the space $L$.

Proof. It is necessary to check that the system of locally defined neighborhoods satisfies the neighborhood axioms NB1 – NB3 from (Kuratowski, 1966).

Let us consider an arbitrary element $l$ from $L$

**NB1:** By definition, any neighborhood $O_N l$ of the element $l$ contains $l$.

**NB2:** Let us consider two neighborhoods $O_N l$ and $O_M l$ of the element $l$ where $N = \{r_t \,;\, t \in K\}$ and $M = \{p_t \,;\, t \in K\}$. Let us take the set $L = \{u_t = \min \{r_t, p_t\} \,;\, t \in K\}$ and show that the intersection $O_M l \cap O_N l$ is equal to the neighborhood $O_L l$ of $l$. Indeed, if an element $h$ belongs both to $O_N l$ and $O_M l$, then

$$\forall t \in K \,(\, q_t(l - h) < r_t)$$

and

$$\forall t \in K \,(\, q_t(l - h) < p_t)$$

Then

$$\forall t \in K \,(\, q_t(l - h) < u_t)$$

Thus, $h$ belongs to $O_L l$ and consequently, $O_M l \cap O_N l \subseteq O_L l$. At the same time, if $h$ belongs to $O_L l$, then

$$\forall t \in K \,(\, q_t(l - h) < u_t)$$

Thus,

$$\forall t \in K \,(\, q_t(l - h) < r_t)$$

and

$$\forall t \in K \,(\, q_t(l - h) < p_t)$$

It means that $h$ belongs to $O_M l \cap O_N l$ and consequently, $O_M l \cap O_N l = O_L l$. Thus, Axiom NB2 is valid.

**NB3:** Let us consider a neighborhood $O_N l$ of the element $l$ and an element $h$ that belongs to $O_N l$. It means that

$$\forall t \in K\ (q_t(l-h) < r_t)$$

Defining $L = \{u_t = k_t\ ;\ t \in K\}$ where $k_t = r_t - q_t(l-h)$ and taking an element $g$ that belongs to $O_L h$, for any $t \in K$ by the triangle inequality for seminorms, we have

$$q_t(l-g) = q_t(l - h + h - g) \le q_t(l-h) + q_t(h-g) < q_t(l-g) + k_t =$$
$$q_t(l-h) + r_t - q_t(l-h) = r_t$$

It means that $g$ belongs to $O_N l$ and consequently, as $g$ is an arbitrary element from $O_L h$, we have $O_L h \subseteq O_N l$. Thus, Axiom NB3 is also valid.

Lemma is proved.

The topology $\upsilon_{lQ}$ defines *Q-local convergence* for elements from **L**.

**Definition 2.4.** If $l$ is an element from **L** and $r$ is a positive real number, then a set of the form

$$O_r l = \{g \in L\ ;\ \exists k \in \mathbf{R}^{++}\ \forall t \in K\ (q_t(l-g) < r - k)\}$$

is called a *uniform neighborhood* of the element $l$.

**Lemma 2.4.** The system of all uniform neighborhoods determines a topology $\upsilon_{uQ}$ called *uniform Q-topology* in the space **L**.

Proof. It is necessary to check that the system of locally defined neighborhoods satisfies the neighborhood axioms NB1 – NB3 from (Kuratowski, 1966).

Let us consider an arbitrary element $l$ from **L**.

**NB1:** By definition, any neighborhood $O_r l$ of the element $l$ contains $l$.

**NB2:** Let us consider two neighborhoods $O_r l$ and $O_p l$ of the element $l$. Let us take $u = \min\{r, p\}$ and show that the intersection $O_M l \cap O_N l$ is equal to the neighborhood $O_L l$ of $l$. Indeed, if an element $h$ belongs both to $O_N l$ and $O_M l$, then

$$\exists k \in \mathbf{R}^{++}\ \forall t \in K\ (q_t(l-h) < r - k)$$

and

$$\exists m \in \mathbf{R}^{++}\ \forall t \in K\ (q_t(l-h) < p - m)$$

Then

$$\exists n \in \mathbf{R}^{++}\ \forall t \in K\ (q_t(l-h) < u - n)$$

where $n = k$ when $u = r$ and $n = m$ when $u = p$.

Thus, $h$ belongs to $O_u l$ and consequently, $O_r l \cap O_p l \subseteq O_u l$. At the same time, if $h$ belongs to $O_u l$, then

$$\exists\, n \in \mathbf{R}^{++} \; \forall\, t \in K \; ( q_t(l - h) < u - n)$$

Thus,

$$\exists\, k \in \mathbf{R}^{++} \; \forall\, t \in K \; ( q_t(l - h) < r - k)$$

and

$$\exists\, m \in \mathbf{R}^{++} \; \forall\, t \in K \; ( q_t(l - h) < p - m)$$

because we can take $k = m = n$.

It means that $h$ belongs to $O_r l \cap O_p l$ and consequently, $O_r l \cap O_p l = O_u l$. Thus, Axiom NB2 is valid.

**NB3:** Let us consider a neighborhood $O_r l$ of the element $l$ and an element $h$ that belongs to $O_r l$. It means that

$$\exists\, k \in \mathbf{R}^{++} \; \forall\, t \in K \; ( q_t(l - h) < r - k)$$

If $p = (\frac{1}{2})k$ and $g$ is an element that belongs to $O_p h$, then there is $u \in \mathbf{R}^{++}$ such that for any $t \in K$ by the triangle inequality for seminorms, we have

$$q_t(l - g) = q_t(l - h + h - g) \leq q_t(l - h) + q_t(h - g) < (r - k) + (p - u) =$$
$$r - k + (\tfrac{1}{2})k - u = r - (\tfrac{1}{2})k - u = r - ((\tfrac{1}{2})k + u)$$

It means that $g$ belongs to $O_r l$ and consequently, as $g$ is an arbitrary element from $O_p h$, we have $O_p h \subseteq O_r l$. Thus, Axiom NB3 is also valid.

Lemma is proved.

The topology $\upsilon_{uQ}$ defines *Q-global convergence* for elements from $L$.

**Proposition 2.8.** The topology $\upsilon_{uQ}$ is stronger than the topology $\upsilon_{lQ}$.

Indeed, any globally defined neighborhood $O_r\, l$ with respect to $Q$ is an open set in the topology $\upsilon_{lQ}$. Namely,

$$O_r l = O_{N_i} l$$

where $N_i = \{\, r_{t,i} = (1 - (1/i))r \,;\, t \in K \,\}$.

Let us consider two sets $Q$ and $T$ of seminorms in $L$.

**Proposition 2.9.** If $Q \subseteq T$, then the topology $\upsilon_{lT}$ is stronger than the topology $\upsilon_{lQ}$ and the topology $\upsilon_{uT}$ is stronger than the topology $\upsilon_{uQ}$

Indeed, any locally defined neighborhood with respect to $T$ is also a locally defined neighborhood with respect to $Q$.

In addition to topologies in $\boldsymbol{L}$, the set $Q = \{q_t \,; t \in K\}$ of seminorms in $\boldsymbol{L}$ defines two types of topology for the set $\boldsymbol{L}^\omega$ of all sequences with elements from $\boldsymbol{L}$: *local topology* and *global* or *uniform topology*.

**Definition 2.5.** If $f = (f_i)_{i \in \omega}$ is a sequence from $\boldsymbol{L}^\omega$ and $N = \{r_t \,; t \in K\}$ is a set of positive real numbers $r_t$, then a set of the form

$$O_N f = \{ (g_i)_{i \in \omega}; (g_i)_{i \in \omega} \in \boldsymbol{L}^\omega \ \& \ \forall t \in K \ \exists k_t \in \boldsymbol{R}^{++} \ \exists n(t) \in \omega \ \forall i > n(t) \ ( q_t(f_i - g_i) < r_t - k_t)\}$$

is called a *locally defined neighborhood* of the sequence $f = (f_i)_{i \in \omega}$.

**Lemma 2.5.** The system of all locally defined neighborhoods determines a topology $\tau_{lQ}$ called *local Q-topology* in the space $\boldsymbol{L}^\omega$.

Proof. It is necessary to check that the system of locally defined neighborhoods satisfies the neighborhood axioms (Kuratowski, 1966).

Let us consider an arbitrary sequence $f$ from $\boldsymbol{L}^\omega$.

**NB1:** By definition, any neighborhood $O_N f$ of the sequence $f$ contains $f$.

**NB2:** Let us consider two neighborhoods $O_N f$ and $O_M f$ of the sequence $f$ where $N = \{r_t \,; t \in K\}$ and $M = \{p_t \,; t \in K\}$. Let us build the set $L = \{l_t = \min \{r_t, p_t\}; t \in K\}$ and show that the intersection $O_M f \cap O_N f$ is equal to the neighborhood $O_L f$ of $f$. Indeed, if a sequence $g$ belongs both to $O_N f$ and $O_M f$, then

$$\forall t \in K \ \exists k_t \in \boldsymbol{R}^{++} \ \exists n(t) \in \omega \ \forall i > n(t) \ ( q_t(f_i - g_i) < r_t - k_t)$$

and

$$\forall t \in K \ \exists h_t \in \boldsymbol{R}^{++} \ \exists m(t) \in \omega \ \forall i > m(t) \ ( q_t(f_i - g_i) < p_t - h_t)$$

Taking some $t \in K$ and defining $u_t = \max \{l_t + k_t - r_t, l_t + h_t - p_t\}$ and $w(t) = \max \{n(t), m(t)\}$, we assume that $r_t < p_t$ and thus, $l_t = r_t$. Then $u_t = \max \{k_t, l_t + h_t - p_t\}$. If $k_t \geq l_t + h_t - p_t$, then $u_t = k_t$ and $q_t(f_i - g_i) < r_t - k_t = l_t - u_t$. If $k_t < l_t + h_t - p_t$, then $u_t = l_t + h_t - p_t$, $p_t = l_t + h_t - u_t$ and $q_t(f_i - g_i) < p_t - h_t = l_t + h_t - u_t - h_t = l_t - u_t$. As the case $r_t > p_t$ is symmetric, we have

$$\forall i > w(t) \ ( q_t(f_i - g_i) < l_t - u_t)$$

Thus, $g$ belongs to $O_L f$ and consequently, $O_M f \cap O_N f \subseteq O_L f$. At the same time, if $g$ belongs to $O_L f$, then

$$\forall t \in K \, \exists v_t \in \mathbf{R}^+ \, \exists w(t) \in \omega \, \forall i > w(t) \, ( \, q_t(f_i - g_i) < l_t - v_t)$$

Thus,

$$\forall i > w(t) \, ( \, q_t(f_i - g_i) < r_t - v_t) \text{ and } \forall i > w(t) \, ( \, q_t(f_i - g_i) < p_t - v_t)$$

It means that $g$ belongs to $O_M f \cap O_N f$ and consequently, $O_M f \cap O_N f = O_L f$. Thus, Axiom NB2 is valid.

**NB3:** Let us consider a neighborhood $O_N f$ of the sequence $f$ and a sequence $h = (h_i)_{i \in \omega}$ that belongs to $O_N f$. It means that

$$\forall t \in K \, \exists k_t \in \mathbf{R}^{++} \, \exists n(t) \, \forall i > n(t) \, ( \, q_t(f_i - h_i) < r_t - k_t)$$

Defining $L = \{l_t = (1/3)k_t \, ; \, t \in K\}$ and $w(t) = \max \{n(t), m(t)\}$ and taking a sequence $g$ that belongs to $O_L h$, we have

$$\forall t \in K \, \exists c_t \in \mathbf{R}^+ \, \exists m(t) \, \forall i > m(t) \, ( \, q_t(h_i - g_i) < l_t - c_t)$$

At the same time, if $w(t) = \max \{n(t), m(t)\}$, then by the triangle inequality (Axiom N3), for any $i > w(t)$, we have

$$q_t(f_i - g_i) = q_t(f_i - h_i + h_i - g_i) \leq q_t(f_i - h_i) + q_t(h_i - g_i) <$$
$$r_t - k_t + (1/3)k_t - c_t = r_t - (2/3)k_t - c_t = r_t - u_t$$

where $u_t = (2/3) k_t + c_t$. It means that $g$ belongs to $O_N f$ and consequently, as $g$ is an arbitrary element from $O_L h$, we have $O_L h \subseteq O_N f$. Thus, Axiom NB3 is valid for the system of locally defined neighborhoods.

Lemma is proved.

The topology $\tau_{lQ}$ defines *Q-local convergence* for sequences of elements from $L$. For instance, when $L$ is the set of all real functions $F(\mathbf{R}, \mathbf{R})$ and $Q = Q_{pt}$, we have the conventional *pointwise convergence* of sequences of real functions (cf. Example 2.3).

Let us assume that the equivalence relation $\approx_Q$ is not trivial in $L^\omega$, i.e., there are, at least, two sequences of elements from $L$ that are equivalent with respect to $\approx_Q$.

**Proposition 2.10.** The topology $\tau_{lQ}$ does not satisfy the axiom $\mathbf{T_0}$ in $L^\omega$, and in particular, $L^\omega$ is not a Hausdorff space with respect to this topology.

Proof. Let us consider two arbitrary sequences $f = (f_i)_{i \in \omega}$ and $g = (g_i)_{i \in \omega}$ of elements from $L$ such that $f \approx_Q g$. Then any neighborhood $O_N f$ of the sequence $f$ contains $g$ and any neighborhood $O_M g$ of the sequence $g$ contains $f$.

Proposition is proved.

**Corollary 2.1.** The space $F(R)^{\omega}$ of all sequences of real functions with the topology of the compact convergence is not a $T_0$-space.

**Definition 2.6.** If $f = (f_i)_{i \in \omega}$ is a sequence from $L^{\omega}$ and $r$ is a positive real number, then a set of the form

$$O_r f = \{(g_i)_{i \in \omega}; (g_i)_{i \in \omega} \in L^{\omega}\ \&\ \exists\, k \in R^{++}\ \forall t \in K\ \exists\, n(t) \in \omega\ \forall i > n(t)\ (q_t(f_i - g_i) < r - k)\}$$

is called a *uniform neighborhood* of the sequence $f = (f_i)_{i \in \omega}$.

When $N = \{r_t = r$ for all $t \in K\}$, Lemma 2.5 implies the following result.

**Lemma 2.6.** The system of all uniform neighborhoods determines a topology $\tau_{uQ}$ called *uniform Q-topology* in the space $L^{\omega}$.

Proposition 2.10 implies the following result.

**Corollary 2.2.** The uniform topology $\tau_{uQ}$ does not satisfy the axiom $T_0$ in $L^{\omega}$, and in particular, $L^{\omega}$ is not a Hausdorff space with respect to this topology.

Proposition 2.9 implies the following result.

**Corollary 2.3.** If $Q \subseteq T$, then the topology $\tau_{lT}$ is stronger than the topology $\tau_{lQ}$ and the topology $\tau_{uT}$ is stronger than the topology $\tau_{uQ}$

The topology $\tau_{uQ}$ defines *Q-uniform convergence* for sequences of elements from $L$.

**Proposition 2.11.** The topology $\tau_{uQ}$ is stronger than the topology $\tau_{lQ}$.

Taking a a set $Q = \{q_t; t \in K\}$ of seminorms in $L$ and the set $L^{\omega}$ of all sequences of elements from $L$, we define two types of topology for the hyperspace $L_{\omega Q}$: local topology and global or uniform topology.

**Definition 2.7.** If $F = Hs_Q(f_i)_{i \in \omega}$ is an element from the hyperspace $L_{\omega Q}$ and $N = \{r_t; t \in K\}$ is a set of real numbers $r_t$, then a set of the form

$$O_N F = \{G;\ \text{there is}\ g = (g_i)_{i \in \omega} \in L^{\omega}\ \text{such that}$$

$$G = \mathrm{Hs}_Q(g_i)_{i \in \omega} \text{ and } \boldsymbol{g} = (g_i)_{i \in \omega} \in \mathrm{O}_N \boldsymbol{f} \text{ where } \boldsymbol{f} = (f_i)_{i \in \omega}\}$$

is called a locally defined neighborhood of the element $F$.

**Lemma 2.7.** The system of all locally defined neighborhoods determines a topology $\delta_{lQ}$ called *local Q-topology* in the set $\boldsymbol{L}_{\omega Q}$.

<u>Proof</u>. Let us consider an arbitrary element $F = \mathrm{Hs}_Q(f_i)_{i \in \omega}$ from the hyperspace $\boldsymbol{L}_{\omega Q}$ and its locally defined neighborhood $\mathrm{O}_N F$ built from the neighborhood $\mathrm{O}_N \boldsymbol{f}$. At first, let us check that the neighborhood $\mathrm{O}_N F$ does not depend on the choice of the representative $(f_i)_{i \in \omega}$ of the element $F$. Taking another representative $\boldsymbol{h} = (h_i)_{i \in \omega}$ of the element $F$ and the locally defined neighborhood $\mathrm{O}^{\circ}_N F$ that is built from the neighborhood $\mathrm{O}_N \boldsymbol{h}$, we show that any element $G$ from $\mathrm{O}_N F$ also belongs to $\mathrm{O}^{\circ}_N F$.

Indeed, if $G$ belongs to $\mathrm{O}_N F$, then there is $(g_i)_{i \in \omega} \in \boldsymbol{L}^{\omega}$ such that $G = \mathrm{Hs}_Q(g_i)_{i \in \omega}$ and $(g_i)_{i \in \omega} \in \mathrm{O}_N \boldsymbol{f}$ where $\boldsymbol{f} = (f_i)_{i \in \omega}$. It means that

$$\forall t \in K \, \exists k_t \in \boldsymbol{R}^{++} \, \exists n(t) \in \omega \, \forall i > n(t) \, ( \, q_t(f_i - g_i) < r_t - k_t)$$

As $F = \mathrm{Hs}_Q(h_i)_{i \in \omega}$, we have $\lim_{i \to \infty} q_t(f_i - h_i) = 0$ for any $t \in K$, and consequently, for some $l_t < k_t$,

$$\exists \, m(t) \, \forall i > m(t) \, ( \, q_t(f_i - h_i) < l_t)$$

By the properties of seminorms,

$$q_t(h_i - g_i) \leq q_t(h_i - f_i) + q_t(f_i - g_i) < r_t - k_t + l_t = r_t - (k_t - l_t) = r_t - p_t$$

where $p_t = k_t - l_t$. Thus, for $c(t) = \max\{n(t), m(t)\}$, we have

$$\forall i > c(t) \, ( \, q_t(f_i - g_i) < r_t - p_t)$$

As a result, we obtain

$$\forall t \in K \, \exists p_t \in \boldsymbol{R}^{++} \, \exists c(t) \in \omega \, \forall i > c(t) \, ( \, q_t(h_i - g_i) < r_t - p_t)$$

It means that the element $G$ belongs to $\mathrm{O}^{\circ}_N F$. Thus, we have proved that $\mathrm{O}_N F \subseteq \mathrm{O}^{\circ}_N F$. As the equivalence relation $\approx_Q$ is symmetric, $\mathrm{O}_N F \subseteq \mathrm{O}^{\circ}_N F$ and consequently, $\mathrm{O}_N F = \mathrm{O}^{\circ}_N F$, i.e., the neighborhood $\mathrm{O}_N F$ does not depend on the choice of the representative of the element $F$.

To prove that the system of all locally defined neighborhoods determines a topology $\delta_{lQ}$ in the set $\boldsymbol{L}_{\omega Q}$, it is necessary to check that the system of locally defined neighborhoods satisfies the neighborhood axioms (cf. Appendix and (Kuratowski 1966)).

**NB1:** By definition, any neighborhood $\mathrm{O}_N F$ of an element $F$ contains $F$.

**NB2:** Let us consider two neighborhoods $O_N\, F$ and $O_M\, F$ of the element $F$ where $N = \{r_t\,;\, t \in K\}$ and $M = \{p_t\,;\, t \in K\}$.

If $\pi_Q\colon L^\omega \to L_{\omega Q}$ is the natural projection, i.e., $\pi_Q(f = (f_i)_{i\in\omega}) = \mathrm{Hs}_Q(f_i)_{i\in\omega}$, then the neighborhood $O_N\, F$ is equal to the projection $\pi_Q(O_N f)$ and the neighborhood $O_M\, F$ is equal to the projection $\pi_Q(O_M f)$. By Lemma 2.5, there is a neighborhood $O_L f$ that belongs to the intersection of both neighborhoods $O_N f$ and $O_M f$. Then its projection $\pi_Q(O_L f)$ is a neighborhood of $F$ that belongs to the intersection of both neighborhoods $O_N\, F = \pi_Q(O_N f)$ and $O_M\, F = \pi_Q(O_M f)$. Thus, Axiom NB2 is valid.

**NB3:** Let us consider a neighborhood $O_N\, F$ of the extrafunction $F$ and an extrafunction $G$ that belongs to $O_N\, F$. It means that there is a sequence $g = (g_i)_{i\in\omega} \in L^\omega$ such that $G = \mathrm{Hs}_Q(g_i)_{i\in\omega}$ and $(g_i)_{i\in\omega} \in O_N f$ where $f = (f_i)_{i\in\omega}$. By Lemma 2.3, there is a neighborhood $O_L g$ such that $O_L g \subseteq O_N f$. Then its projection $\pi_Q(O_L g)$ is a neighborhood of $G$ and $O_L\, G = \pi_Q(O_L g) \subseteq \pi_Q(O_N f) = O_N\, F$. Thus, Axiom NB3 is also valid for the system of locally defined neighborhoods.

Lemma is proved.

**Theorem 2.3.** The topology $\delta_{lQ}$ satisfies the axiom $T_2$, and thus, $L_{\omega Q}$ is a Hausdorff space with respect to this topology.

<u>Proof.</u> Let us consider two arbitrary elements $F = \mathrm{Hs}_Q(f_i)_{i\in\omega}$ and $G = \mathrm{Hs}_Q(g_i)_{i\in\omega}$ from $L_{\omega Q}$. If $F \neq G$ in $L_{\omega Q}$, then any sequences $f = (f_i)_{i\in\omega} \in F$ and $g = (g_i)_{i\in\omega} \in G$ satisfy the following condition: there is a seminorm $q$ from $Q$ and a positive number $k$ such that for any $n \in \omega$, there is an $i > n$ such that $q(f_i - g_i) > k$. This condition makes it possible to choose an infinite set $M$ of natural numbers such that for any $m \in M$ the inequality $q(f_i - g_i) > k$ is valid.

Let us take $l = k/4$ and consider two uniform neighborhoods $O_l f$ and $O_l g$ of the sequences $f$ and $g$ in $L^\omega$. If $\pi_Q\colon L^\omega \to L_{\omega Q}$ is the natural projection, then the projections $\pi_Q(O_l f)$ and $\pi_Q(O_l g)$ of these neighborhoods $O_l f$ and $O_l g$ are neighborhoods $O_l\, F$ and $O_l\, G$ of the elements $F$ and $G$ with respect to the topology $\delta_{lQ}$. By construction, $\pi_Q(O_l f) \cap \pi_Q(O_l g) = \varnothing$. To prove this, we suppose that this is not true. Then there is a an element $H = \mathrm{Hs}_Q(h_i)_{i\in\omega}$ from $L_{\omega Q}$ that is an element of the set $\pi_Q(O_l f) \cap \pi_Q(O_l g)$. This implies that there are sequences $u = (u_i)_{i\in\omega}$ and $v = (v_i)_{i\in\omega}$ from $L^\omega$ for which $\pi_Q(u) = \pi_Q(v) = H$, $u \in O_l f$ and $v \in O_l g$.

The equality $\pi_Q(u) = \pi_Q(v)$ implies that for the chosen seminorm $q$ and number $k$, the following condition (*) is valid: $\exists m \in \omega\, \forall i > m\, (q(u_i - v_i) < k/3)$. The set $M$, which is determined above, is infinite. So, there is $j \in M$ such that it is greater than $m$ and $q(u_j - v_j) \geq q(f_j - g_j) - q(f_j - u_j) - q(g_j - v_j) \geq k - k/4 - k/4 = k/2 > k/3$ because $f_j - g_j = f_j - u_j + u_j - v_j - g_j + v_j$ and by properties of seminorms, $q(f_j - g_j) \leq q(f_j - u_j) + q(u_j - v_j) + q(g_j - v_j)$, $q(f_j - u_j) < l - r < l = k/4$, $q(g_j - v_j) < l - p < l = k/4$, and $q(v_j - g_j) = q(g_j - v_j)$. This contradicts the condition (*), according to which $q(u_j - v_j) < k/3$.

Consequently, the assumption is not true, and $\pi_Q(O_l f) \cap \pi_Q(O_l g) = O_l F \cap O_l G = \emptyset$.

Theorem is proved because $F = \mathrm{Hs}_Q(f_i)_{i \in \omega}$ and $G = \mathrm{Hs}_Q(g_i)_{i \in \omega}$ are arbitrary elements from $L_{\omega Q}$.

The topology $\delta_{lQ}$ defines *Q-local convergence* of elements from $L_{\omega Q}$.

**Definition 2.8.** If $F = \mathrm{Hs}_Q(f_i)_{i \in \omega}$ is an element from the hyperspace $L_{\omega Q}$ and $r$ is a real number, then a set of the form

$$O_r F = \{G; \text{ there is } \mathbf{g} = (g_i)_{i \in \omega} \in L^\omega \text{ such that }$$
$$G = \mathrm{Hs}_Q(g_i)_{i \in \omega} \text{ and } \mathbf{g} = (g_i)_{i \in \omega} \in O_r \mathbf{f} \text{ where } \mathbf{f} = (f_i)_{i \in \omega}\}$$

is called a *uniform neighborhood* of the element $F$.

When $N = \{r_t = r \text{ for all } t \in K\}$, Lemma 2.7 implies the following result.

**Lemma 2.8.** The system of all uniform neighborhoods determines a topology $\delta_{uQ}$ called *uniform Q-topology* in the hyperspace $L_{\omega Q}$.

The topology $\delta_{uQ}$ defines *Q-uniform convergence* of elements from the hyperspace $L_{\omega Q}$.

Theorem 2.3 implies the following result.

**Corollary 2.4.** The uniform topology $\delta_{uQ}$ satisfies the axiom $T_2$, and thus, the hyperspace $L_{\omega Q}$ is a Hausdorff space with respect to this topology.

Proposition 2.9 implies the following result.

**Corollary 2.5.** If $Q \subseteq T$, then the topology $\delta_{lT}$ is stronger than the topology $\delta_{lQ}$ and the topology $\delta_{uT}$ is stronger than the topology $\delta_{uQ}$

Proposition 2.11 implies the following result.

**Proposition 2.12.** The topology $\delta_{uQ}$ is stronger than the topology $\delta_{lQ}$.

## 3. Bundles with a Hyperspace Base

Taking a (topological) vector space $L$, its hyperspace $L_{\omega Q}$ and the natural projection $\pi_Q\colon L^\omega \to L_{\omega Q}$, we obtain the *fiber bundle* $\mathbf{BD}L_{\omega Q} = (L^\omega, \pi_Q, L_{\omega Q})$ with the base $L_{\omega Q}$, bundle space $L^\omega$ and the fiber $F = \{\, f \in L^\omega;\, f \approx_Q \mathbf{0} = (0, 0, 0, \ldots)\,\}$ as . $\mathbf{BD}L_{\omega Q}$ is a *vector bundle*.

We will use this bundle for differentiation in the hyperspace $L_{\omega Q}$.

**Example 3.1**. Taking the space $R$ of all real numbers, their absolute values as the norm and applying the construction of the hyperspace, we obtain the real-number vector bundle $\mathbf{BD}R_\omega = (R^\omega, \pi_R, R_\omega)$.

**Example 3.2**. Taking the space $C$ of all complex numbers, their absolute values as the norm and applying the construction of the hyperspace, we obtain the complex-number vector bundle $\mathbf{BD}C_\omega = (C^\omega, \pi_C, C_\omega)$.

**Example 3.3**. Taking the space $F(R, R)$ of all real functions and the hyperspace $F(R, R)_{\omega Q\mathrm{pt}}$ or $F(R, R_\omega)$ of all restricted real pointwise extrafunctions (cf. Example 2.3), we obtain the restricted real pointwise extrafunction vector bundle $\mathbf{BD}F(R, R)_{\omega Q\mathrm{pt}} = (F(R, R)^\omega, \pi_{R\mathrm{pt}}, F(R, R)_{\omega Q\mathrm{pt}})$.

**Example 3.4**. Taking the space $F(R, C)$ of all functions from $R$ to $C$ and the hyperspace $F(R, C)_{\omega Q\mathrm{pt}}$ or $F(R, C_\omega)$ of all restricted complex valued pointwise extrafunctions (cf. Example 2.4), we obtain the restricted real-complex pointwise extrafunction vector bundle $\mathbf{BD}F(R,C)_{\omega Q\mathrm{pt}} = (F(R, C)^\omega, \pi_{RC\mathrm{pt}}, F(R, C)_{\omega Q\mathrm{pt}})$.

**Example 3.5**. Taking the space $F(C, C)$ of all complex functions and the hyperspace $F(C, C)_{\omega Q\mathrm{pt}}$ of all restricted complex pointwise extrafunctions (cf. Example 2.5), we obtain the restricted complex pointwise extrafunction vector bundle $\mathbf{BD}F(C, C)_{\omega Q\mathrm{pt}} = (F(C, C)^\omega, \pi_{C\mathrm{pt}}, F(C, C)_{\omega Q\mathrm{pt}})$.

**Example 3.6**. Taking a set $\mathbf{K}$ of test functions, the space $F(R, R)$ of all real functions and the hyperspace $F(R, R)_{\omega Q\mathbf{K}}$ of real extended distributions (cf. Example 2.6), we obtain the real extended distributions vector bundle $\mathbf{BD}F(R, R)_{\omega Q\mathbf{K}} = (F(R, R)^\omega, \pi_{R\mathbf{K}}, F(R, R)_{\omega Q\mathbf{K}})$.

**Example 3.7**. Taking the space $C(R, R)$ of all continuous real functions and the hyperspace $C(R, R)_{\omega Q\mathrm{comp}}$ or $\mathrm{Comp}(R, R_\omega)$ of all real compactwise extrafunctions (cf. Example 2.7), we obtain the real compactwise extrafunction vector bundle $\mathbf{BD}C(\mathbf{R}, \mathbf{R})_{\omega Q\text{comp}} = (C(\mathbf{R}, \mathbf{R})^{\omega}, \pi_{R\text{comp}}, C(\mathbf{R}, \mathbf{R})_{\omega Q\text{comp}})$.

If $\mathbf{B} = (E, p, B)$ is a fiber bundle, then a *section* in $B$ is a mapping $r: B \to E$ such $\pi_Q \circ r$ is the identity mapping of $B$.

If $\mathbf{B} = (E, p, B)$ is a vector bundle, then the set SK($\mathbf{B}$) of all sections of the bundle $\mathbf{B}$ is a vector space. Indeed, the fiber $F$ is a vector space and for each element $x$ from $B$, it is possible to define operations in the vector space $p^{-1}(b) = F$. So, we have *addition* of sections $r + q$ is defined as

$$(r + q)(x) = r(x) + q(x) \text{ for any } x \in B;$$

and *multiplication* of a section $r$ by an element $a$ from the field $ar$ is defined as

$$(ar)(x) = r(ax) \text{ for any } x \in B \text{ and } a \in F.$$

In addition, we consider several types of vector bundle sections.

**Definition 3.1.** A section $r: B \to E$ is called:

a) *additive* if $r(x + y) = r(x) + r(y)$ for any $x, y \in B$;

b) *uniform* if $r(ax) = a(r(x))$ for any $x \in B$ and $a \in F$;

c) *linear* if it is additive and uniform;

d) *continuous* if $r$ is a continuous mapping;

e) stable if $r(\pi_Q(f)) = f$ for any stable element $f$ from $L^{\omega}$.

Not all sections are additive, uniform, linear or continuous.

**Example 3.8**. Taking the real-number vector bundle $\mathbf{BDR}_{\omega} = (\mathbf{R}^{\omega}, \pi_R, \mathbf{R}_{\omega})$, we consider the section $r: \mathbf{R}^{\omega} \to \mathbf{R}_{\omega}$ where for any stable hypernumber $\alpha = \text{Hn}(a, a, a, \ldots, a, \ldots)$, we have

$$r(\alpha) = (a, a, a, \ldots, a, \ldots) \text{ when } a \neq 2$$

and

$$r(2, 2, 2, \ldots, 2, \ldots) = (2 + (1), 2 + (½), 2 + (1/3), \ldots, 2 + (1/n), \ldots)$$

This section is not additive because $r(1) + r(1) \neq r(2)$ and is not uniform because $2r(1) \neq r(2)$.

At the same time, there is a general method to build linear sections. Taking a fiber bundle $\mathbf{BDL}_{\omega Q} = (L^{\omega}, \pi_Q, L_{\omega Q})$ and a basis $V$ of the vector space $L_{\omega Q}$, we obtain the subbundle $E_V =$ $(E, \pi_V, V)$ where $E = \pi_Q^{-1}(V)$ and $\pi_D$ is the restriction of $\pi_Q$ onto the set $E$. Usually subbundles with the same base as the initial bundle are considered. However, here we use a more general definition of a subbundle, which corresponds to the category of bundles with different bases.

**Proposition 3.1.** a) Any section of the bundle $E_V$ has a unique extension to a linear section of the bundle $\mathbf{BDL}_{\omega Q}$.

b) The restriction of a linear section of the bundle $\mathbf{BDL}_{\omega Q}$ is a section of the bundle $E_V$.

c) A section of the bundle $\mathbf{BDL}_{\omega Q}$ is linear if and only if it is an extension of a section of the bundle $E_V$.

Proof. (a) Let us consider a section $q: V \to E$. in the bundle $E_V = (E, \pi_V, V)$. If $F \in L_{\omega Q}$, then $F = \sum_{i=1}^{n} a_i H_i$ where $a_i \in F$ and $H_i \in V$ ($i = 1, 2, 3, \ldots, n$). This allows us to define

$$r(F) = \sum_{i=1}^{n} a_i\, q(H_i)$$

This definition is correct because the decomposition $F = \sum_{i=1}^{n} a_i H_i$ is unique. As $F$ is an arbitrary element from $L_{\omega Q}$, we obtain a section of the bundle $\mathbf{BDL}_{\omega Q}$. By construction, this section is linear. (a) is proved.

(b) is true because $E = \pi_Q^{-1}(V)$.

(c) follows from (a) and (b).

Proposition is proved.

Taking an element $a$ from the field $F$, we define the operator $mt_a$ in $L$ by the formula:

$$mt_a(x) = ax$$

Taking an element $e$ from the space $L$, we define the operator $ad_e$ in $L$ by the formula:

$$ad_e(x) = x + e$$

There are natural extensions $mt_a{}^*$ of the operator $mt_a$ and $ad_e{}^*$ of the operator $ad_e$ to the spaces $L^\omega$ and $L_{\omega Q}$:

$$mt_a{}^*\,(\mathrm{Hs}_Q(f_i)_{i \in \omega}) = \mathrm{Hs}_Q(af_i)_{i \in \omega}$$

and

$$ad_e{}^*\,(\mathrm{Hs}_Q(f_i)_{i \in \omega}) = \mathrm{Hs}_Q(f_i + e)_{i \in \omega}$$

as well as

$$mt_a{}^*(f_i)_{i \in \omega} = (af_i)_{i \in \omega}$$

and

$$ad_e^* \, (f_i)_{i \in \omega} = (f_i + e)_{i \in \omega}$$

**Proposition 3.2.** If $r$ is a section of the bundle $\mathbf{BD}L_{\omega Q}$, then:

a) $mt_a \circ r \circ mt_{a^{-1}}$ is also a section of the bundle $\mathbf{BD}L_{\omega Q}$.

b) $ad_e \circ r \circ ad_{-e}$ is also a section of the bundle $\mathbf{BD}L_{\omega Q}$.

This shows that the multiplicative group of $F$ and the additive group of $L$ act on the set SK($\mathbf{BD}L_{\omega Q}$) of all sections of the bundle $\mathbf{BD}L_{\omega Q}$.

To define differentiation in the hyperspace $L_{\omega Q}$, we assume that $L$ is a *differential vector space* over a field $F$ with a partial multiplication, i.e., there is a linear mapping (differential operator)

$$\partial : L \to L$$

that satisfies the Leibniz product law when the product $f \cdot g$ is defined, i.e., if $f, g \in L$, then

$$\partial(f \cdot g) = (\partial f) \cdot g + f \cdot (\partial g)$$

When multiplication is defined for all elements from $L$, the space $L$ is called a *differential algebra* (Kaplansky, 1957; Kolchin, 1973).

In many interesting cases, there is a *standard injection* $\mu_F : F \to L$, i.e., an injection such that $af = \mu_F(a) \cdot f = f \cdot \mu_F(a)$ for any element $a$ from $F$ and any element $f$ from $L$. For instance, when is the space of all real functions, the standard injection of the field $R$ or real numbers assign constant functions to real numbers. Note that it is possible to multiply any element from $L$ by any element from the standard injection and $\mu_F(0)$ is equal to zero $\mathbf{0}$ of the vector space $L$.

**Proposition 3.3.** The derivative $\partial(\mu_F(a))$ of any element $\mu_F(a)$ is equal to $\mathbf{0}$.

Proof. By definition and the Leibniz product law,

$$\partial(\mu_F(a)) = \partial(1 \cdot \mu_F(a)) = \partial(\mu_F(1) \cdot \mu_F(a)) =$$

$$(\partial(\mu_F(1)))\mu_F(a) + \mu_F(1)(\partial(\mu_F(a))) = (\partial(\mu_F(1)))\mu_F(a) + \partial(\mu_F(a))$$

Thus,

$$(\partial(\mu_F(1)))\mu_F(a) + \partial(\mu_F(a)) = \partial(\mu_F(a))$$

and consequently,

$$(\partial(\mu_F(1)))\mu_F(a) = \mathbf{0}$$

Then as $\mu_F$ is a homomorphism,

$$\partial(\mu_F(1)) = \partial(\mu_F(a)) \cdot 1 = (\partial(\mu_F(1)))\mu_F(a \cdot a^{-1}) = (\partial(\mu_F(1)))\mu_F(a) \, \mu_F(a^{-1}) = \mathbf{0}$$

So,

$$\partial(\mu_F(a)) = \partial(\mu_F(a \cdot 1)) = \partial(a \cdot \mu_F(1)) = a \cdot \partial(\mu_F(1)) = \mathbf{0}$$

Proposition is proved.

In a natural way, the operator $\partial$ is extended to the space $\boldsymbol{L}^{\omega}$ when we apply it separately to each coordinate, namely, for $\boldsymbol{f} = (f_i)_{i \in \omega}$, we define

$$\partial \boldsymbol{f} = (\partial f_i)_{i \in \omega}$$

**Example 3.9.** If $\boldsymbol{f} = (x^n)_{n \in \omega}$, then $\partial \boldsymbol{f} = \mathbf{d/d}x\, \boldsymbol{f} = (d/dx\, x^n)_{n \in \omega} = (nx^{n-1})_{n \in \omega}$.

As $\partial$ satisfies the identity

$$\partial(af + bg) = a(\partial f) + b(\partial g)$$

where $a, b \in \boldsymbol{F}$ and $f, g \in \boldsymbol{L}$. Consequently, we have a similar identity in $\boldsymbol{L}^{\omega}$:

$$\partial(a\boldsymbol{f} + b\boldsymbol{g}) = a(\partial \boldsymbol{f}) + b(\partial \boldsymbol{g})$$

where $\boldsymbol{f}, \boldsymbol{g} \in \boldsymbol{L}^{\omega}$, i.e., $\partial$ is a differential operator in $\boldsymbol{L}^{\omega}$.

Proposition 3.3 implies the following example.

**Proposition 3.4.** The derivative $\partial \boldsymbol{a}$ of any element $\boldsymbol{a}$ from $\mu_F(\boldsymbol{F})^{\omega}$ is equal to $\mathbf{0}$.

**Corollary 3.1.** The derivative $\partial \boldsymbol{f}$ of any element $\boldsymbol{f} = (f_n)_{n \in \omega}$ from $C(\boldsymbol{R}, \boldsymbol{R})^{\omega}$ is equal to $\mathbf{0}$ if all $f_n$ are constant functions.

We can try to construct a similar extension of the differential operator $\partial$ for the space $\boldsymbol{L}_{\omega Q}$. However, in a general case, it is impossible to do this in a regular way as the following examples demonstrate.

**Example 3.10.** Let us consider the hyperspace $C(\boldsymbol{R}, \boldsymbol{R})_{\omega Q_{\text{comp}}}$ (or $\text{Comp}(\boldsymbol{R}, \boldsymbol{R}_{\omega})$) of all real compactwise extrafunctions (cf. Example 2.7 and (Burgin, 2004)) and the sequence of functions $f_1(x) = (½)\sin 2x, \ldots, f_n(x) = (½)^n \sin 2^n x, \ldots$. Then the compactwise extrafunction $\text{Ec}(f_n)_{n \in \omega}$ is identically equal to zero, i.e., $\text{Ec}(f_n)_{n \in \omega} = \text{Ec}(g_n)_{n \in \omega}$ where $g_n(x) = 0$ for all $n$ and all $x$. Consequently, $\text{Ec}(d/_{dx}g_n)_{n \in \omega} = \text{Ec}(g_n)_{n \in \omega}$ is also the compactwise extrafunction identically equal to zero. At the same time, $d/_{dx}f_n = \sin 2^n x$ for all $n = 1, 2, 3, \ldots$. Thus, $\text{Ec}(d/_{dx}f_n)_{n \in \omega} \neq \text{Ec}(d/_{dx}g_n)_{n \in \omega}$. More exactly, $\text{Ec}(d/_{dx}f_n)_{n \in \omega}$ is the compactwise extrafunction identically equal to one.

**Example 3.11.** Let us consider the hyperspace $F(\boldsymbol{R}, \boldsymbol{R})_{\omega Q_{\text{pt}}}$ (or $F(\boldsymbol{R}, \boldsymbol{R}_{\omega})$) of all real pointwise extrafunctions (cf. Example 2.3 and (Burgin, 2004)) and the sequence of functions $f_1(x) = \sin x, f_2(x) = (½)\sin 2x, \ldots, f_n(x) = (½)^n \sin 2^n x, \ldots$. Then the pointwise extrafunction $\text{Ep}(f_n)_{n \in \omega}$ is identically equal to zero, i.e., $\text{Ep}(f_n)_{n \in \omega} = \text{Ep}(g_n)_{n \in \omega}$ where $g_n(x) = 0$ for all $n$ and all $x$.

Consequently, $\text{Ep}(d/d_x g_n)_{n\in\omega} = \text{Ep}(g_n)_{n\in\omega}$ is also the pointwise extrafunction identically equal to zero. At the same time, $d/d_x f_n = \sin 2^n x$ for all $n = 1, 2, 3, \ldots$. Thus, $\text{Ep}(d/d_x f_n)_{n\in\omega} \neq \text{Ep}(d/d_x g_n)_{n\in\omega}$. More exactly, $\text{Ep}(d/d_x f_n)_{n\in\omega}$ is the pointwise extrafunction identically equal to one.

Let us consider the fiber bundle $\mathbf{BD}L_{\omega Q} = (L^{\omega}, \pi_Q, L_{\omega Q})$ where $L$ is a differential vector space over a field $F$ and its section $r: L_{\omega Q} \to L^{\omega}$.

**Definition 3.2.** The differential operator

$$\partial/\partial r : L_{\omega Q} \to L_{\omega Q}$$

is defined by the following formula

$$\partial/\partial_r \text{Hs}_Q(f_i)_{i\in\omega} = \pi_Q(\partial(r(\text{Hs}_Q(f_i)_{i\in\omega})))$$

The element $\partial/\partial_r \text{Hs}_Q(f_i)_{i\in\omega}$ from the hyperspace $L_{\omega Q}$ is called a *sectional derivative* of the element $\text{Hs}_Q(f_i)_{i\in\omega}$.

**Example 3.12**. If $F \in L_{\omega Q}$ and $r(F) = (x^n)_{n\in\omega}$, then $\partial/\partial_r F = \text{Hs}_Q(nx^{n-1})_{n\in\omega}$.

Proposition 3.4 implies the following example.

**Proposition 3.5.** The derivative $\partial/\partial_r F$ of a element $F$ from $L_{\omega Q}$ is equal to $\mathbf{0}$ if $r(F)$ belongs to $\mu_F(F)^{\omega}$.

We remind that an extrafunction is called constant if it is equal to a hypernumber.

**Corollary 3.2.** The derivative $\partial/\partial_r F$ of an extrafunction $F$ is equal to $\mathbf{0}$ if $r(F)$ is a hypernumber.

**Remark 3.1.** Examples 3.10 and 3.11 show that an element from the hyperspace $L_{\omega Q}$ can have different sectional derivatives, i.e., in a general case, a sectional derivative of an element from a hyperspace is not unique.

**Theorem 3.1.** If $r$ is an additive section of the bundle $\mathbf{BD}L_{\omega Q}$, then for any elements $F$ and $G$ from the hyperspace $L_{\omega Q}$, we have

$$\partial/\partial_r (F + G) = \partial/\partial_r F + \partial/\partial_r G$$

<u>Proof</u>. Let us take two elements $F = \text{Hs}_Q(f_i)_{i\in\omega}$ and $G = \text{Hs}_Q(g_i)_{i\in\omega}$ from the hyperspace $L_{\omega Q}$. Then by definition

$$\partial/\partial_r (F + G) = \pi_Q(\partial(r(F + G))) = \pi_Q(\partial(r(F) + r(G)))$$

as the section $r$ is additive.

$$\pi_Q(\partial(r(F) + r(G))) = \pi_Q(\partial(r(F)) + \partial(r(G)))$$

as $\partial$ is a linear mapping. Then by Theorem 2.2, we have

$$\pi_Q(\partial(r(F) + \partial(r(G))) = \pi_Q(\partial(r(F))) + \pi_Q(\partial(r(G)))$$

Consequently,

$$\partial/\partial r\, (F + G) = \partial/\partial r\, F + \partial/\partial r\, G$$

Theorem is proved.

**Corollary 3.3.** If $r$ is a additive section, then for any element $F$ from the hyperspace $L_{\omega Q}$, we have

$$\partial/\partial r\, (ad_e*(F)) = ad_{\partial e}*\, (\partial/\partial r F)$$

where $ad_e*$ is a natural extension of the operator $ad_e$ to the spaces $L^\omega$ and $L_{\omega Q}$.

**Theorem 3.2.** If $r$ is a uniform section of the bundle $\mathbf{BD}L_{\omega Q}$, then for any element $f$ from the hyperspace $L_{\omega Q}$ and any real number $c$, we have

$$\partial/\partial r\, (c\cdot F) = c\cdot(\partial/\partial r\, F)$$

Proof. By definition

$$\partial/\partial r\, (c\cdot F) = \pi_Q(\partial(r(c\cdot F))) = \pi_Q(\partial(cr(F)))$$

as the section $r$ is uniform.

$$\pi_Q(\partial(cr(F))) = \pi_Q(\partial(c\cdot \partial(r(F))))$$

as $\partial$ is a linear mapping. Then by Theorem 2.2, we have

$$\pi_Q(c\cdot \partial(r(F))) = c\cdot \pi_Q(\partial(r(F)))$$

Consequently,

$$\partial/\partial r\, (c\cdot F) = c\cdot(\partial/\partial r\, F)$$

Theorem is proved.

**Corollary 3.4.** If $r$ is a linear section, then for any elements $F$ and $G$ from the hyperspace $L_{\omega Q}$ and any real numbers $c$ and $d$, we have

$$\partial/\partial r\, (c\cdot F + d\cdot G) = c\cdot (\partial/\partial r\, F) + d\cdot (\partial/\partial r\, G)$$

**Corollary 3.5.** If $r$ is a uniform section, then for any element $F$ from the hyperspace $L_{\omega Q}$, we have

$$\partial/\partial r\, (mt_a*(F)) = mt_a*\, (\partial/\partial r F)$$

where $mt_a{}^*$ is a natural extension of the operator $mt_a$ to the spaces $L^\omega$ and $L_{\omega Q}$.

**Proposition 3.6.** If $r$ is a stable section and $F$ is a stable element from the hyperspace $L_{\omega Q}$, then $\partial/\partial_r F$ is also a stable element.

However, in many interesting cases, $L$ is not always a differential vector space and the differential operator $\partial$ in $L$ is partial. Then the domain $D$ of the operator $\partial$ is a subspace of $L$ as the operator $\partial$ is linear. We also assume that there is a subspace $H$ of $L$ elements of which are continuous in some sense. For instance, when $L$ is a space of functions (mappings) from a topological space into another topological space, then $H$ consists of continuous functions (mappings). We have another example when continuous elements are fuzzy continuous functions (mappings) in the sense of (Burgin, 2008).

However, having a partial differential operator $\partial$ in $L$, we can extend it to the partial differential operator $\partial$ in $L^\omega$. Namely, if $f = (f_i)_{i \in \omega} \in L^\omega$, then

$$\partial f = \begin{cases} (g_i)_{i \in \omega} \text{ where } g_i = 0 \text{ for all } i = 1, 2, 3, \ldots, n-1 \\ \text{and } g_i = \partial f_i \text{ for all } i = n, n+1, n+2, n+3, \ldots \\ \text{when are all } f_i \text{ differentiable for } i = n, n+1, n+2, n+3, \ldots \\ \\ \text{undefined otherwise} \end{cases}$$

Let us consider the fiber bundle $\mathbf{BD}L_{\omega Q} = (L^\omega, \pi_Q, L_{\omega Q})$, its section $r: L_{\omega Q} \to L^\omega$ and $F \in L_{\omega Q}$. Then

$$\partial/\partial_r F = \pi_Q(\partial(r(F)))$$

**Lemma 3.1.** An element $F \in L_{\omega Q}$ is differentiable if it has a differentiable representation, i.e., $F = \pi_Q(f)$ where $f$ is differentiable in $L^\omega$.

Let us consider the fiber bundle $\mathbf{BD}L_{\omega Q} = (L^\omega, \pi_Q, L_{\omega Q})$ and its section $r: L_{\omega Q} \to L^\omega$.

**Definition 3.3.** A section $r: L_{\omega Q} \to L^\omega$ is called:

a) *differentially represented* if all elements in any sequence from the space $r(L_{\omega Q})$ are differentiable, i.e., they belong to $D$;

b) *continuously represented* if all elements in any sequence from the space $r(L_{\omega Q})$ are continuous, i.e., they belong to $H$.

Lemma 3.1 implies the following result, which shows importance of differentially represented sections.

**Corollary 3.6.** If the fiber bundle $\mathbf{BD}L_{\omega Q}$ has a differentially represented section, then all elements from $L_{\omega Q}$ are differentiable.

**Proposition 3.7.** For any compactwise extrafunction $F = \mathrm{Ec}(f_n)_{n \in \omega}$ and any differentially represented section $r$ of the bundle $C(\mathbf{R}, \mathbf{R})^{\omega}Q_{\mathrm{comp}}$, if $\partial/\partial_r F = \mathbf{0}$ where $\mathbf{0}$ is the identically equal to zero compactwise extrafunction, then $F = \boldsymbol{a}$ where $\boldsymbol{a}$ is the compactwise extrafunction that is identically equal to the number $a$.

Proof. By definition, if $r(\mathrm{Hs}_Q(f_i)_{i \in \omega}) = (h_i)_{i \in \omega}$, then

$$\partial/\partial_r F = \partial/\partial_r \mathrm{Ec}_Q(f_i)_{i \in \omega} = \pi_Q(\mathbf{d/d}_x(r(\mathrm{Ec}_Q(f_i)_{i \in \omega}))) = \pi_Q(\mathbf{d/d}_x(h_i)_{i \in \omega}) = \pi_Q((\mathrm{d/d}_x h_i)_{i \in \omega}) = \mathbf{0}$$

By the definition of real compactwise extrafunctions (Burgin, 2004), it means that the sequence $(\mathrm{d/d}_x h_i)_{i \in \omega}$ uniformly converges to zero in each interval.

Let us consider two arbitrary real numbers $c$ and $d$ with $d < c$. By definition, $F(c) = \mathrm{Ec}(f_n(c))_{n \in \omega} = \mathrm{Ec}(h_n(c))_{n \in \omega}$ and $F(d) = \mathrm{Ec}(f_n(d))_{n \in \omega} = \mathrm{Ec}(h_n(d))_{n \in \omega}$. As each function $h_i$ is differentiable, by the Mean Value Theorem (cf., for example, (Burgin, 2008)), there is a point $e_i$ in the interval $[c, d]$ such that

$$(h_i(c) - h_i(d))/(c - d) = h_i'(e_i) = \mathrm{d/d}_x h_i(e_i)$$

Thus, $\lim_{i \to \infty}(h_i(c) - h_i(d)) = 0$ because $\lim_{i \to \infty} \mathrm{d/d}_x h_i(e_i) = 0$ and $c - d$ is a constant. By the definition of hypernumbers,

$$F(c) = \mathrm{Ec}(h_n(c))_{n \in \omega} = \mathrm{Ec}(f_n(d))_{n \in \omega} = F(d)$$

As $c$ and $d$ are arbitrary real numbers, $F$ is a constant compactwise extrafunction.

Proposition is proved.

**Remark 3.2.** Examples 3.10 and 3.11 show that in contrast to ordinary differentiable functions, a constant compactwise extrafunction can have non-zero sectional derivatives.

**Proposition 3.8.** For any pointwise extrafunction $F = \mathrm{Ec}(f_n)_{n \in \omega}$ and any differentially represented section $r$ of the bundle $\mathbf{BD}F(\mathbf{R}, \mathbf{R})_{\omega}Q_{\mathrm{pt}}$, if $\partial/\partial_r F = \mathbf{0}$ where $\mathbf{0}$ is the identically equal to zero pointwise extrafunction, then $F = \boldsymbol{a}$ where $\boldsymbol{a}$ is the pointwise extrafunction that is identically equal to the number $a$.

Proof is similar to the proof of Proposition 3.7.

**Remark 3.3.** Example 3.10 shows that in contrast to ordinary differentiable functions, a constant pointwise extrafunction can have non-zero sectional derivatives.

The Weierstrass theorem allows us to prove the following result.

**Theorem 3.3.** The bundle $\mathbf{BD}C(\mathbf{R}, \mathbf{R})_{\omega Q\text{comp}}$ of compactwise extrafunctions has differentially represented sections.

This theorem has interesting corollaries.

**Corollary 3.7.** All compactwise extrafunctions are differentiable.

Proof is based on Theorem 3.3 and Corollary 3.1.

Thus, even if differentiation in $L$ is partial, it is possible that differentiation in $L_{\omega Q}$ is total. This explains why extrafunctions are so useful for differential equations, allowing one to solve equations that do not have solutions even in distributions.

**Definition 3.4.** An element $H$ from the hyperspace $F(\mathbf{R}, \mathbf{R})_{\omega Q\text{pt}}$ is called:

a) *differentially represented* if there is a sequence $f = (f_i)_{i \in \omega}$ such that $H = \mathrm{Hs}_{Q\text{pt}}(f_i)_{i \in \omega}$ and all functions in $f$ are differentiable;

b) *continuously represented* if there is a sequence $f = (f_i)_{i \in \omega}$ such that $H = \mathrm{Hs}_{Q\text{pt}}(f_i)_{i \in \omega}$ and all functions in $f$ are continuous.

As the sum of two continuous functions is a continuous function and multiplication of a continuous function gives a continuous function, we have the following result.

**Proposition 3.9.** The set $F^C(\mathbf{R}, \mathbf{R})_{\omega Q\text{pt}}$ of all continuously represented pointwise extrafunctions is a subspace of the hyperspace $F(\mathbf{R}, \mathbf{R})_{\omega Q\text{pt}}$.

Proposition 3.5 allows us to build the bundle $\mathbf{BD}F^C(\mathbf{R}, \mathbf{R})_{\omega Q\text{pt}} = (C(\mathbf{R}, \mathbf{R})^{\omega}, \upsilon_{R\text{pt}}, F^C(\mathbf{R}, \mathbf{R})_{\omega Q\text{pt}})$, where $C(\mathbf{R}, \mathbf{R})^{\omega}$ is the space of all sequences of continuous real functions, as the subbundle of the bundle $\mathbf{BD}F(\mathbf{R}, \mathbf{R})_{\omega Q\text{pt}}$.

Then the Weierstrass theorem allows us to prove the following result.

**Theorem 3.4.** The bundle $\mathbf{BD}F^C(\mathbf{R}, \mathbf{R})_{\omega Q\text{pt}}$ of continuously represented pointwise extrafunctions has differentially represented sections.

**Corollary 3.8.** All compactwise extrafunctions are differentiable.

Proof is based on Theorem 3.4 and Corollary 3.1.

In a similar way, we can prove the following results.

**Theorem 3.5.** The bundle **BD***D*(***R***) of extended distributions has differentially represented sections.

**Corollary 3.9.** All extended distributions are differentiable.

Proof is based on Theorem 3.5 and Corollary 3.1.